\documentclass{conm-p-l}

\usepackage{graphicx}
\usepackage{epsf}
\usepackage{latexsym}

\title[Cayley Trick and Products of Simplices]
{The Cayley Trick and
 Triangulations of \\
Products of Simplices
}

\author{Francisco Santos}
 \address{
  Universidad de Cantabria,
  Departamento de Matem\'aticas,
  Estad\'{\i}stica y Computa\-ci\'on, Facultad de Ciencias,
  39005 Santander, SPAIN.
 }

\thanks{
This paper was finished in the fall of 2003, while I was a member
of the Mathematical Sciences Research Institute (Berkeley, CA), 
 supported by MSRI and the Spanish Ministry of Education.
My research is also supported by grant BFM2001-1153 of the Spanish
Ministry of Science and Technology.
}

\date{March 5, 2004}

\theoremstyle{plain}
\newtheorem{theoremintro}{Theorem}
\newtheorem{theorem}{Theorem}[section]
\newtheorem{lemma}[theorem]{Lemma}
\newtheorem{proposition}[theorem]{Proposition}
\newtheorem{corollary}[theorem]{Corollary}

\theoremstyle{definition}
\newtheorem{definition}[theorem]{Definition}

\theoremstyle{remark}
\newtheorem{remark}[theorem]{Remark}
\newtheorem{question}[theorem]{Question}

\newcommand{\tconv}{\operatorname{tconv}}

\newcommand{\vertices}{\operatorname{vertices}}
\newcommand{\height}{\operatorname{height}}

\newcommand{\CP}{{\mathbb C\mathbb P}}
\newcommand{\TP}{{\mathbb T \mathbb P}}

\newcommand{\reals}{{\mathbb R}}
\newcommand{\R}{{\mathbb R}}
\newcommand{\naturals}{{\mathbb N}}

\renewcommand{\H}{{\mathcal H}}
\newcommand{\C}{{\mathcal C}}

\begin{document}

\email{ santosf@unican.es}
\urladdr{ http://personales.unican.es/santosf}

\begin{abstract}
We use the Cayley Trick to study polyhedral subdivisions of the
product of two simplices. For arbitrary (fixed) $l\ge 2$,
we show that the numbers of regular and non-regular triangulations
of $\Delta^l\times\Delta^k$
grow, respectively, as $k^{\Theta(k)}$ and $2^{\Omega(k^2)}$.

For the special case of $\Delta^2\times \Delta^k$, we relate
triangulations to certain class of lozenge tilings. This allows us to
compute the exact number of triangulations up to $k=15$, show that the
number grows as $e^{\beta k^2/2 + o(k^2)}$ where $\beta\simeq
0.32309594$ and prove that the set of all triangulations is
connected under geometric bistellar flips. The latter has as a
corollary that the toric Hilbert scheme of the determinantal ideal of
$2\times 2$ minors of a $3\times k$ matrix is connected, for every
$k$.

We include ``Cayley Trick pictures'' of all the triangulations of
$\Delta^2\times \Delta^2$ and $\Delta^2\times \Delta^3$, as well as
one non-regular triangulation of $\Delta^2\times \Delta^5$ and another of
$\Delta^3\times \Delta^3$.
\end{abstract}

\maketitle

\section*{Introduction}
\label{sec:intro}

The {\it polyhedral Cayley Trick}\, gives a canonical bijection
between {\em mixed subdivisions} of the Minkowski sum $\sum_{i=1}^k P_i$
of several polytopes $P_1,\dots,P_k\in \R^d$ and all polyhedral
subdivisions of a certain polytope $\C(P_1,\dots,P_k)\subset
\R^d\times \R^{k-1}$ called the {\it Cayley embedding} of
$P_1,\dots,P_k$. The correspondence was first developed by
Sturmfels~\cite{Sturmfels94} for the case of coherent subdivisions, and
then generalized to all subdivisions by Huber et al.~\cite{HuRaSa}.

Originally, the trick was devised as a way of understanding and
computing fine (i.e., minimal with respect to refinement) mixed
subdivisions, taking advantage of 
the much deeper knowledge and specific software
that exists for triangulations. But the trick can be also used in reverse, to 
understand triangulations of the $d+k-1$-dimensional polytope
$\C(P_1,\dots,P_k)$ in terms of a $d$-dimensional object. 
This is what we do here.

Specially interesting is the case when all the $P_i$'s are copies of a
simplex $\Delta^{l-1}$. Then, the Cayley Trick relates polyhedral
subdivisions of $\Delta^{k-1}\times \Delta^{l-1}$ to mixed subdivisions
of the dilation $k\Delta^{l-1}$. If, moreover, we fix $l=3$, then
the mixed subdivisions we have to study are essentially the same 
as lozenge tilings of a triangle of size $k$.
Using this interpretation we prove:

\begin{theoremintro}
\label{thm:main}
\begin{enumerate}
\item The graph of flips between triangulations of
$\Delta^2\times \Delta^{k}$ is
connected (Theorem~\ref{thm:connected})  
and it has diameter $\Theta(k^2)$ (Corollary~\ref{coro:diameter}).

\item  (Theorem~\ref{thm:asymptotics}) The number of triangulations
of $\Delta^2\times \Delta^k$ grows as
\[
e^{\frac{3}{\pi}L\left(\frac{\pi}{3}\right) k^2 + o(k^2)},
\]
where $L(x)$ is the Lobachevsky function
\[
L(x)=-\int_{0}^{x} \log |2\sin t| dt.
\]

\item The number of triangulations of $\Delta^2\times\Delta^{k-1}$, for 
$k=1,\dots,16$ is $k!$ times the number shown in Table~\ref{table:numbers}. 
\end{enumerate}

\begin{table}[htb]
\rm
\hrule
\medskip
\begin{center}
\begin{tabular}{|cc|}
\hline
$k$ &          lozenge tilings of $k\Delta^2$ \cr
\hline
1 &                          1 \cr
2 &                          3 \cr
3 &                         18 \cr
4 &                        187 \cr
5 &                       3135 \cr
6 &                      81462 \cr
7 &                    3198404 \cr
8 &                  186498819 \cr
\hline
\end{tabular}
\hskip .5 cm
\begin{tabular}{|cc|}
\hline
$k$ &              lozenge tilings of $k\Delta^2$ \cr
\hline
9 &                15952438877 \cr
10&              1983341709785 \cr
11&            355891356876534 \cr
12&          91655826195854811 \cr
13&       33726014269095727260 \cr
14&    17665249123640876125464 \cr
15& 13130399067698641838496272 \cr
16&13813411778618644581617635925\cr
\hline
\end{tabular}
\hfill
\end{center}
\caption{Tilings of $k\Delta^2$ into $k$ triangles and ${k\choose 2}$
lozenges. Times $k!$, this is the number of
triangulations of $\Delta^2\times\Delta^{k-1}$}
\label{table:numbers}
\hrule
\end{table}

\end{theoremintro}

The number $\frac{3}{\pi}L\left(\frac{\pi}{3}\right)\simeq 0.323$ 
that appears in part (2) of this statement is
the maximum asymptotic normalized entropy of lozenge tilings of a
planar region, as computed in~\cite{CoKePr}.
The exact number of triangulations of $\Delta^2\times\Delta^{k}$
that appears in part (3) had previously been computed only up to
$k=5$~\cite{Loera,Rambau-TOPCOM}.
Part (1) of the statement
is interesting for two reasons. On the one hand, there are not many 
examples where the graph of flips is known to be connected. Essentially, only
the case of dimension at most 2 (classical), codimension at 
most 3~\cite{AzaSan} and that of cyclic polytopes~\cite{Rambau-cyclic}.
On the other hand, since all triangulations of products of simplices are 
unimodular, the graph of flips has a very direct interpretation in
toric algebraic geometry; see Theorem~\ref{thm:torichilbert}
below.

Besides the results in Theorem~\ref{thm:main}, the Cayley Trick allows to picture triangulations of $\Delta^2\times
\Delta^k$ as $2$-dimensional objects. We include pictures of all
(non-isomorphic) triangulations of $\Delta^2\times \Delta^2$ and 
$\Delta^2\times \Delta^3$
(Figures~\ref{fig:twobytwo} and~\ref{fig:35classes}). Also, of
non-regular triangulations of $\Delta^2\times \Delta^5$ and
$\Delta^3\times \Delta^3$ and
of a non-regular coarse subdivision of $\Delta^2\times \Delta^7$
(Figures~\ref{fig:nonregular} 
and~\ref{fig:3x3nonregular}). $\Delta^2\times \Delta^5$ and $\Delta^3\times
\Delta^3$ are the minimal products of simplices that have non-regular
triangulations. To the best of our knowledge, our subdivision 
of $\Delta^2\times \Delta^7$ is the first known non-regular coarse
subdivision of a product of simplices.

\medskip

Subdivisions and triangulations of $\Delta^{k-1}\times
\Delta^{l-1}$ are interesting from several perspectives. 
They have been studied for their
own sake in~\cite{BabBil,Loera}and~\cite[Sect. 7.3.D]{GKZbook} and
they have been used as
building blocks to find efficient triangulations of high dimensional
cubes~\cite{Haiman,OrdSan} or to find disconnected flip-graphs~\cite{Santos-noflips,Santos-toric}.
Also, $\Delta^{k-1}\times \Delta^{l-1}$ is an example of a {\em totally
unimodular} polytope; that is,
a lattice polytope all the simplices of which have
the same volume. This implies it is equidecomposable~\cite{Bayer},
i.e., all its triangulations have the same $f$-vector. From the
$f$-vector of any of its triangulations one can recover, for example,
its Erhart polynomial.

From a more algebraic point of view, the toric ideal associated
to the vertex set of  $\Delta^{k-1}\times \Delta^{l-1}$ is a
fundamental object. It is the determinantal ideal generated by the $2\times
2$-minors of a $k\times l$-matrix, and the variety associated to it is
the {\em Segre embedding} of $\CP^{k-1}\times \CP^{l-1}$ in $\CP^{kl-1}$. The
study of this ideal has connections to
enumeration, sampling and optimization
for contingency tables and transportation problems
(cf.~\cite[Chapter 5]{Sturmfels-book}). 

In this context, part (1) of Theorem~\ref{thm:main} has the following 
algebro-geometric application:
Every lattice point set defines a toric (binomial)
ideal and a toric Hilbert scheme. Contrary to standard Hilbert
schemes, toric Hilbert schemes are sometimes non-connected ~\cite{Santos-toric}. 
In fact, for a totally unimodular polytope, connectivity
of the corresponding toric Hilbert scheme is equivalent to connectivity of the graph
of triangulations
(cf.~\cite{MacTho} and Theorem 10.13 in~\cite{Sturmfels-book}).
Hence:

\begin{theoremintro}
\label{thm:torichilbert}
The toric Hilbert scheme of the determinantal ideal of $2\times 2$ minors of a
$3\times k$ matrix is connected.
\end{theoremintro}

Another recent source of interest in triangulations of $\Delta^{k-1}\times\Delta^{l-1}$
comes from tropical geometry. 
Develin and Sturmfels~\cite{DevStu} have proved that the combinatorial
types of point configurations with $k$ points in the
tropical space of dimension $l-1$ are in bijective correspondence to
the {\em regular} subdivisions of $\Delta^{k-1}\times \Delta^{l-1}$.
Via this relation, we prove:

\begin{theoremintro}
[Corollary~\ref{coro:nonregular2}]
\label{thm:tropical}
For every fixed $l\ge 2$, the number of regular subdivisions of 
$\Delta^{l}\times\Delta^{k}$ grows as $2^{\Theta(k\log k)}$
while the number of non-regular subdivisions grows as 
$2^{\Omega(k^2)}$.
\end{theoremintro}

The structure of the paper is as follows: After a first section with
preliminaries on the Cayley Trick, Section~\ref{sec:labeling} shows how to study mixed
subdivisions from a purely geometric point of view. The results hold
for arbitrary mixed subdivisions, but have a specially simple form for
dilations of a simplex. Sections~\ref{sec:lozenge1} and~\ref{sec:lozenge2} contain our main results, on triangulations of
$\Delta^2\times \Delta^{k-1}$. The first explores the relation between
mixed subdivisions of $k\Delta^2$ and lozenge tilings, and the second 
uses it to prove Theorem~\ref{thm:main}.
Finally, Section~\ref{sec:tropical} des\-cribes the relation between
regular subdivisions of $\Delta^l\times\Delta^k$ and tropical
polytopes in $l$-space, and proves Theorem~\ref{thm:tropical}.
\medskip

\section{The Cayley Trick}

Except for Theorem~\ref{thm:free}, the contents of this section are
special cases of results from~\cite{HuRaSa}.

\subsection{Polyhedral subdivisions}

Let $P$ be a polytope in $\R^d$. A {\em cell} of $P$ is any sub-polytope of
the same dimension as $P$ and whose vertices are a subset of those of $P$. A
{\em polyhedral subdivision} of $P$ is any family of cells which cover $P$ and
intersect properly,  meaning that $B\cap B'$ is a common face of
both $B$ and $B'$ for every pair of cells $B$ and $B'$ in the subdivision.
This definition is a special case of the definition of polyhedral subdivisions of a point
configuration~\cite{Triangbook,Loera,GKZbook,Ziegler}.

Polyhedral subdivisions form a poset under the refinement relation 
\[
S \le S' \qquad \Longleftrightarrow \qquad 
\forall B\in S \quad \exists B'\in S' \ : \ B\subseteq B'.
\]
The minimal elements in this poset, 
i.e., those subdivisions whose cells are all simplices,
are the {\em triangulations} of $P$.
The unique maximal element in the poset is the \emph{trivial
subdivision}, which has the whole polytope $P$ as its only cell.
The subdivisions which refine only
the trivial one are called {\em coarse}.

A polyhedral subdivision is called \emph{regular} (or, sometimes,
coherent) if it can be obtained as the orthogonal projection of the
lower facets of a $d+1$-dimensional polytope. Equivalently,
$S$ is regular if there 
is a height function $h:\operatorname{vertices}(P) \to \R$ such that 
for every cell $B\in S$ the points $\{(v,h(v)):v\in B\}\subset
\R^{d+1}$ lie in a hyperplane that passes strictly below  all other
points $\{(w,h(w)) : w\not\in B\}$.

\subsection{Mixed subdivisions}
Let $P_1,\ldots,P_k\,\,\subset\R^d$ be convex polytopes.
We in principle do not need to assume them to be full-dimensional, but
we require their Minkowski sum to be. This Minkowski sum is defined as:
\begin{displaymath}
\left\{ \,x_1+\cdots+x_k : x_i\in P_i 
\right\}\mbox{.}
\end{displaymath}

\begin{definition}
  A {\em Minkowski cell} of the Minkowski sum $\sum_{i=1}^kP_i$ is any
  full-dimensional polytope $B=\sum_{i=1}^k B_i$, where each $B_i$ is a
  (perhaps not full-dimensio\-nal) polytope with vertices among those of $P_i$.
  A {\em mixed subdivision} of $\sum_{i=1}^kP_i$ is any family 
$S$ of Minkowski
  cells which cover $\sum P_i$ and intersect properly
  as Minkowski sums, meaning that for any
  two cells $B=\sum B_i$ and $B'=\sum B'_i$ in $S$ and for every
  $i\in\{1,\dots,k\}$, the polytopes $B_i$ and $B'_i$ intersect properly
(their intersection is a face of both).
\end{definition}

\begin{remark}
We have intentionally decided to slightly abuse notation in these definitions 
and in the rest of the paper, to simplify it 
and to make the geometry more apparent.
Indeed, although we speak of  ``mixed subdivisions of $P=\sum P_i$'', the concept depends 
not only on the polytope $P$ but also on the particular Minkowski decomposition
of it that we are given. The same occurs in the 
definition of pro\-per intersection of $B$ and $B'$. To resolve this ambiguity,
{\em every time we write a Minkowski sum, the expression
$\sum_{i=1}^k Q_i$ should be formally understood as an 
ordered $k$-tuple of polytopes $(Q_1,\dots,Q_k)$}. Only in sentences like
``a family of Minkowski cells covers $\sum P_i$'' we are referring to the underlying polytope
resulting from the sum.

In particular, if there are different ways
of obtaining a certain subpolytope of $P$ as a sum of subpolytopes from the
$P_i$'s, we consider these as different Minkowski cells and we have to specify
which of them we are using in a particular mixed subdivision. For example,
part (a) of 
Figure~\ref{fig:twosquares} shows a mixed subdivision of the Minkowski sum of
two equal squares. Below the figure, each of the six Minkowski cells is
expressed as a sum of sub-polytopes of the squares $a_1 a_2 a_3 a_4$ and
$b_1 b_2 b_3 b_4$.
Clearly, the exchange of all the $a$'s and $b$'s in these
expressions would still provide a (different) mixed subdivision, with the same
picture. On the other hand, considering $B_4 = a_4 + 
b_1 b_3 b_4$ and $B_5 = a_1 a_2 a_3 + b_4$ would not provide a mixed
subdivision, since these two Minkowski cells do not intersect properly. 

\begin{figure}[htb]
\begin{center}
    \includegraphics[height=4.2cm]{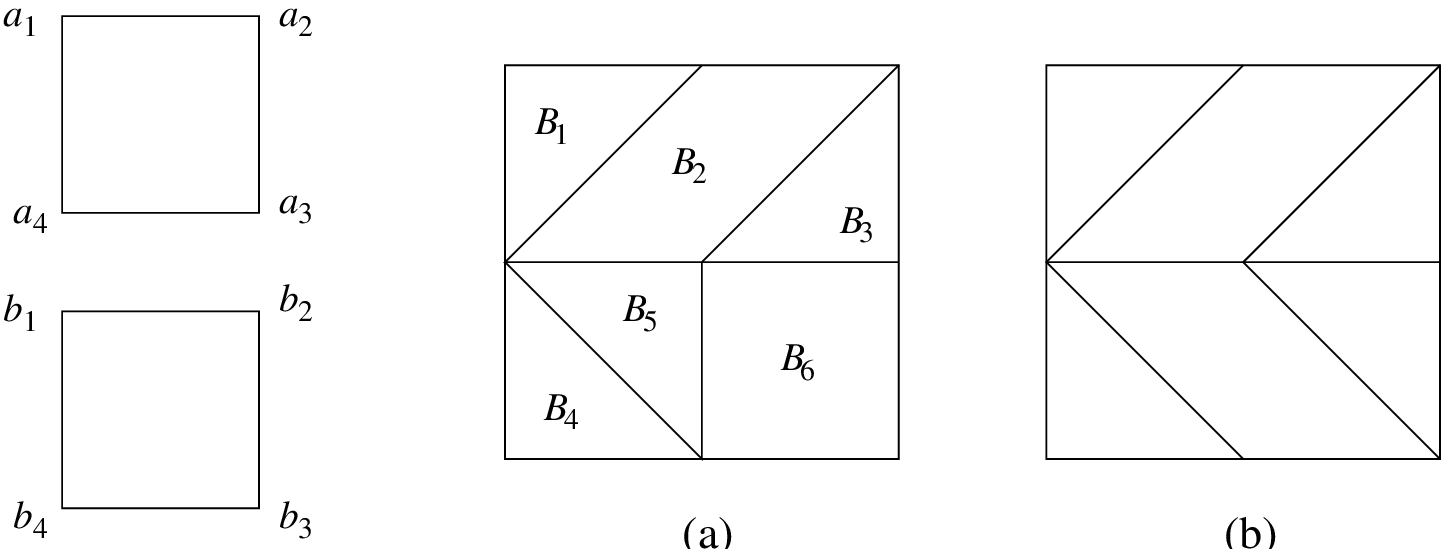}
\small
\[
 B_1 = a_1 a_2 a_4 + b_1, \quad
 B_2 = a_2 a_4 + b_1 b_2, \quad
 B_3 = a_2 a_3 a_4 + b_2, 
\]
$
 B_4 = a_4 + b_1 b_3 b_4, \quad
 B_5 = a_4 + b_1 b_2 b_3, \quad
 B_6 = a_3 a_4 + b_2 b_3. 
$
  \caption{(a) A mixed subdivision of the Minkowski sum of two equal squares. 
  (b) A subdivision which is not mixed}
\label{fig:twosquares}
\end{center}
\end{figure}

Part (b) of the figure shows another reason why the ``labeling'' of cells as Minkowski sums is important.
There, the same Minkowski sum of two equal squares is
decomposed into cells which intersect properly as polytopes and which
individually can be considered Minkowski cells, but which cannot be labeled 
in a way that makes them intersect properly in the
Minkowski sense. In other words, the picture is not compatible with any mixed subdivision.
\end{remark}

\begin{remark}
Our treatment of mixed subdivision is formally different, but equivalent to,
the one in~\cite{HuRaSa}. There, in order to keep track of the Minkowski decompositions
of cells, instead of speaking of the sum of polytopes
$\sum P_i$ the authors speak of the sum of point sets $\sum A_i$, where each $A_i$
is the vertex set of $P_i$. $\sum A_i$ is considered as a multiset, where different 
``points'' (different Minkowski sums) may have identical geometric coordinates.

For the reader familiar with the theory of fiber polytopes (cf. Chapter 10 in~\cite{Ziegler}),
mixed subdivisions of $\sum P_i$ can also be described as the subdivisions which are $\Pi$-induced, 
where $\Pi$ is the natural projection 
\[
\Pi: \Delta^{\# \vertices(P_1)-1}\times \cdots\times \Delta^{\# \vertices(P_k)-1} \to P_1 + \cdots + P_k.
\]
This is Lemma 2.4 in~\cite{HuRaSa}.
\end{remark}


\medskip

As in the case of subdivisions, mixed subdivisions form a poset. 
A Minkowski cell $B=\sum B_i$ is smaller than (or contained in) 
another one $B'=\sum B'_i$ if for every $i$ we have $B_i\subseteq B'_i$. We
write $B \le B'$ if this happens. This induces the following
refinement relation among mixed subdivisions:
\[
S \le S' \qquad \Longleftrightarrow \qquad
\forall B\in S \quad \exists B'\in S' \ : \ B\le B'.
\]

The unique maximal element in this poset is again the trivial
subdivision with only one Minkowski cell $B=\sum P_i$. The minimal
elements are called {\em fine mixed subdivisions}. In Proposition~\ref{prop:fine}
we will see that fine mixed subdivisions are characterized in terms of
what Minkowski cells they use. This is analogous to the fact that ``fine
polyhedral subdivisions'' (i.e., triangulations)
of a polytope are the polyhedral subdivisions whose cells are simplices.
As in the case of
subdivisions, we call {\em coarse mixed subdivisions} those which only
refine the trivial one. They do not have an easy intrinsic
characterization, as far as we know.

A mixed subdivision of $P_1,\dots,P_k$ is called \emph{coherent} if
it can be obtained as the orthogonal projection of the
lower facets of a $d+1$-dimensional Minkowski sum 
$\tilde P_1 +\cdots+\tilde P_k\subset\reals^{d+1}$, where each 
$\tilde P_i$ orthogonally projects to the corresponding $P_i$.
Equivalently, $S$ is regular if there 
are height functions $h_i:\operatorname{vertices}(P_i) \to \R$ such that 
for every cell $B=\sum B_i\in S$ the points 
\[
\{(v_1,h_1(v_1))+\cdots + (v_k,h_k(v_k)) :
              v_i\in B_i\text{ for all } i\}\subset
\R^{d+1}
\]
lie in a hyperplane that passes strictly below  all other points 
\[
\{(w_1,h_1(w_1))+\cdots + (w_k,h_k(w_k)) :
              w_i\not\in B_i \text{ for some } i\}.
\]

%
%
%
%
%

\subsection{The Cayley Trick}
As before, let $P_1,\dots,P_k$ be polytopes in $\reals^d$. Let
$e_1,\dots,e_k$ be an affine basis in $\R^{k-1}$ and let
$\mu_i:\R^d\to\R^d\times\R^{k-1}$ be the inclusion 
$\mu_i(x)=(x,e_i)$. We define
\begin{displaymath}
  \C(P_1,\dots,P_k):=\hbox{conv}(\cup_{i=1}^k \mu_i(P_i)).
\end{displaymath}
We call $\C(P_1,\dots, P_k)$ the {\em Cayley embedding} of $P_1$,\dots,$P_k$.
We will primarily be interested in the Cayley
embedding of $k$ copies of the same polytope $P$. Clearly, it  equals
$P\times \Delta^{k-1}$, where $\Delta^{k-1}$ is the simplex with vertices
$e_1,\dots,e_k$.

The Cayley Trick (see~\cite{HuRaSa} for a full exposition) is a poset
isomorphism between polyhedral subdivisions of the Cayley embedding and mixed
subdivisions of the Minkowski sum. More precisely, observe that for any choice
of affine coordinates $(\lambda_1,\dots,\lambda_k)$ in the simplex
$\Delta^{k-1}$, the intersection of $\C(P_1,\dots,P_k)$ with the
affine subspace $\R^d\times\{\sum \lambda_i e_i\}$ equals the ``weighted
Minkowski sum'' $\sum \lambda_i P_i$. Choosing
$\lambda_1=\cdots=\lambda_k=1/k$ gives a (scaled) copy of the standard
Minkowski sum. What the Cayley Trick says is that:

\begin{theorem}[Cayley Trick~{\cite[Theorem 3.1]{HuRaSa}}]
\ 
\vskip -\baselineskip
\begin{enumerate}
\item For any polyhedral subdivision of $\C(P_1,\dots,P_k)$, intersecting its
  cells with $\R^d\times\{\sum \frac{1}{k} e_i\}$ we get a mixed subdivision of
  $\frac{1}{k}\sum P_i$.
\item This correspondence is a poset isomorphism between polyhedral subdivisions
  of $\C(P_1,\dots,P_k)$ and mixed subdivisions of $\sum P_i$, and
bijects regular subdivisions of  $\C(P_1,\dots,P_k)$ to coherent mixed
subdivisions of $\sum P_i$.
\end{enumerate}
\end{theorem}

Part (1), and hence one direction of part (2), are straightforward:
every full-dimensional cell in a subdivision of $\C(P_1,\dots,P_k)$ is
itself a Cayley embedding  $\C(B_1,\dots,B_k)$ of certain subpolytopes 
$B_i\subseteq P_i$, and hence its intersection with  $\R^d\times\{\sum
\frac{1}{k} e_i\}$ is a Minkowski cell in $\frac{1}{k}\sum P_i$. 
The other direction, that every mixed subdivision arises in this way,
can be easily proved using different values of
$\lambda=(\lambda_1,\dots,\lambda_k)\in \Delta^{k-1}$: The ``Minkowski
intersection property'' of the cells in a mixed subdivision $S$ of 
$P_1+\cdots+P_k$ guarantees that $S$ induces a mixed subdivision of
$\sum \lambda_i P_i$ for every $\lambda$ (by just replacing each
Minkowski cell $B=\sum B_i$ by its weighted version
$\lambda B:=\sum \lambda_i B_i$) and the cells
$\C(B_1,\dots,B_k)=\cup_{\lambda \in \Delta^{k-1}}\lambda B$ 
for $B\in S$ form the desired
polyhedral subdivision of $\C(P_1,\dots,P_k)$.

\subsection{The case $P_1=\cdots=P_k$}
When all the polytopes $P_i$ are copies of
a single polytope $P\subset \reals^d$, then
$\C(P,\dots,P)=P\times\Delta^{k-1}$.
Hence, the Cayley Trick is a
bijection between subdivisions of the product $P\times \Delta^{k-1}$
and mixed subdivisions of the dilation $kP$ of $P$.  The
main topic of this paper is to take advantage of this fact in order to
study triangulations of $P\times \Delta^{k-1}$, via fine mixed
subdivisions of $P+\cdots +P$.
Here comes a first result in this direction. 

\begin{theorem}
In the action of the affine symmetry group of $ \Delta^{k-1}$ on
triangulations of $P\times\Delta^{k-1}$, every orbit has 
exactly $k!$ elements (that is, the action is free).
In particular, the number of triangulations of  $P\times\Delta^{k-1}$
is divisible by $k!$.
\label{thm:free}
\end{theorem}

\begin{proof}
The action of the permutation group, when regarded in the
corresponding mixed subdivision of $kP$, amounts to a mere permutation
of the Minkowski summands of every Minkowski cell, without affecting
the cell itself (as a polytope). In particular, if some permutation
sends a mixed subdivision $S$ to itself, it will send every Minkowski
cell of $S$ to itself.  Hence, it suffices to prove that for every 
fine mixed subdivision $S$
and every permutation $\sigma$ there is
a Minkowski cell $B=\sum B_i$ that is not invariant under reordering of
the summands (that is to say, in which $B_i\ne B_{\sigma(i)}$ for some
$i$).

Finding such cells is easy. Take any $i$ such that $i\ne \sigma(i)$ 
and let $B=\sum B_i$ be a cell in which $B_i$ is
full-dimensional. These are easy to find in the Cayley embedding: they
are the ones with a full-dimensional (in $P$) intersection with 
the face $P\times \{e_i\}$ of $P\times\Delta^{k-1}$, where $e_i$ is
the $i$th vertex of  $\Delta^{k-1}$.
\end{proof}

Theorem~\ref{thm:free} can be proved directly (and easily) in the world
of triangulations of $P\times \Delta^{k-1}$. What we want to emphasize
is that the proof is much more transparent using the Cayley Trick. In
particular, we have not found this result in the bibliography about
triangulations of the product of two simplices~\cite{Loera,BabBil}.

\section{The labeling of a mixed subdivision}
\label{sec:labeling}

\subsection{General case}
We said earlier (and illustrated with Figure~\ref{fig:twosquares})
that there may be different ways in which a given polyhedral
subdivision of $P_1 + \cdots + P_k$ can be labeled as a mixed
subdivision. In this section we address the problem of what
additional information is needed to make the mixed subdivision unique.
We start with a straightforward, but useful, observation:

\begin{lemma}
\label{lemma:facets}
Let $I\subset\{1,\dots,k\}$ be a subset of indices. Let $S$ be a mixed
subdivision of $P_1+\cdots+P_k$. Then, the following is a mixed
subdivision of $\sum_{i\in I} P_i$:
\[
S|_{I}:=\{\sum_{i\in I} B_i : \sum_{i\in \{1,\dots,k\}} B_i \in S \}.
\]
\end{lemma}

\begin{proof}
To see that the cells in $S|_{I}$ cover $\sum_{i\in I} P_i$ one can
use a limiting process: for each
$\lambda=(\lambda_1,\dots,\lambda_k)\in (\reals_+)^k$
the cells $\{\lambda B : B\in S\}$ form a mixed subdivision of
$\sum_{i\in \{1,\dots,n\}} \lambda_i P_i$. If we fix $\lambda_i=1$ for
the indices in $S$ and make all other $\lambda$'s go to zero, the
Minkowski sum tends to $\sum_{i\in I} P_i$.
That the cells intersect properly in the Minkowski sense is straightforward.
\end{proof}

This result can be easily understood in the Cayley embedding: 
$\C(P_i : i\in I)$ is a face of $\C(P_1,\dots,P_k)$ and, certainly,
every polyhedral subdivision of $\C(P_1,\dots,P_k)$ induces a
subdivision of it.

\begin{theorem}
\label{thm:labels}
Let $S$ be a polyhedral subdivision of $P_1+\cdots+P_k$ and suppose
that a labeling of it as a mixed subdivision exists. Assume further
that all the $P_i$ have the same dimension. If we know:
\begin{enumerate}
\item The subdivisions $S_i:=S|_{\{i\}}$ induced by $S$ in each of the
individual $P_i$'s, and
\item Which cell of $S$ collapses to each full-dimensional cell of
each $S_i$,
\end{enumerate}
then the whole mixed subdivision can be recovered from that data
(that is to say, $S$ is the only mixed subdivision compatible with
that information and there is an algorithm to recover the whole labeling).
\end{theorem}

\begin{proof}
We work with one $i$ at a time. That is to say, we are going to fix
$i$ and show how to recover the $i$-th summand $B_i$ of a certain cell
$B\in S$. We argue by induction on the codimension of $B_i$. The
hypotheses give us the case of codimension zero.

So, suppose that we already know the $i$-th summand of all cells
for which that summand has dimension $d'+1$. Let $B_i$ be a
$d'$-dimensional cell in $S_i$. Our goal is to determine all maximal
cells of $S$ whose $i$-th component is precisely $B_i$:

-- First, if 
a cell $C$ of $S$ has been determined to have as $i$-th summand a cell
$C_i\in S_i$ that contains $B_i$ (hence as a face), then the
faces of $C$ in any of the directions defined by the
(relatively open) normal cone) of $B_i$ in $C_i$ will have $B_i$ as
their $i$-th component.

-- Second, if a full-dimensional cell $C'$ 
has been determined by the previous rule to have a face with $i$-th
summand equal to $B_i$ then its $i$-th summand contains $B_i$ as a
face. In particular, it either has been determined already or equals $B_i$.

Our claim is that applying these two rules we can recover all the
cells whose $i$-th summand is $B_i$. A way to see it is in the
process that gives $S_i$ as the limit of mixed subdivisions of
$\lambda (P_1+\cdots+P_k) + (1-\lambda) P_i$ then $\lambda$ goes to
zero. When $\lambda$ is very close to zero, all the cells whose $i$-th
summand is $B_i$ become very close to $B_i$ itself and certainly we
can go from any of them to one that becomes close to the
afore-mentioned $C_i$ traversing only cells whose $i$-th summand is $B_i$.
\end{proof}

\subsection{Fine and pure subdivisions}
Lemma~\ref{lemma:facets} makes it
easy to show that fine mixed subdivisions can
be characterized in terms of what Minkowski cells they are allowed to
use:

\begin{proposition}
\label{prop:fine}
A mixed subdivision $S$ is fine if and only if 
for every Min\-kows\-ki cell $B=\sum B_i$ in $S$, 
the $B_i$'s are all
simplices and $\dim(B)=\sum \dim(B_i)$ (in other words, the $B_i$'s lie in
independent affine subspaces).
\end{proposition}

\begin{proof}
The ``if'' direction is straightforward: If all the cells are as
claimed then no mixed cell can be properly contained in one of $S$
and, hence, $S$ is minimal. For the ``only if'' direction, a direct
proof can be given but it is simpler to use the Cayley Trick. The
Cayley Trick implies that fine mixed subdivisions are those for which
the associated polyhedral of $\C(P_1,\dots,P_k)$ is a
triangulation. In other words, those that use only Minkowski cells
$B=\sum B_i$ whose associated Cayley cells
$\C(B_1,\dots,B_k)$ are simplices. This condition is easily seen to be
equivalent to the one in the statement.
\end{proof}

We will refer to the special mixed cells described in 
Proposition~\ref{prop:fine} as {\em fine mixed cells}. 

\begin{corollary}
\label{coro:labels}
For a fine mixed subdivision, the information in (1) and (2) of
Theorem~\ref{thm:labels} can be recovered if we know,
for each $i$, what maximal cells of $S$ have a full-dimensional 
$i$-th summand.
\end{corollary}

\begin{proof}
In a fine mixed subdivision, if a cell has a full-dimensional summand
then all the other summands are 0-dimensional (points). Hence,
the information we are given is what maximal cells of $S$ form the
subdivision $S_i$, for every $i$, except we
are not told how to arrange them to subdivide $P_i$. If we show how to
do that, then Theorem~\ref{thm:labels} gives the rest.
It is actually enough to find out, for every cell $B$ of $S$, what is
(up to translation) the $i$th summand $B_i$ of that cell. If we know
this, we know how to write $B=B_i +C_i$ for each cell (where the $C_i$ is the
sum of all the other summands) and we can scale down the $C_i$
components of all cells to recover (in the limit) the subdivision
$P_i$.

Suppose then that we have identified (up to translation) all the
$i$-th summands of dimension greater than a certain $d'$. 
Exactly as in Theorem~\ref{thm:labels}, it is then possible to propagate
along $S$ all the $d'$ dimensional faces of those summands, hence
getting all the $i$-th summands of dimension $d'$.
\end{proof}

\begin{figure}[htb]
\begin{center}
\includegraphics[height=1.9in]{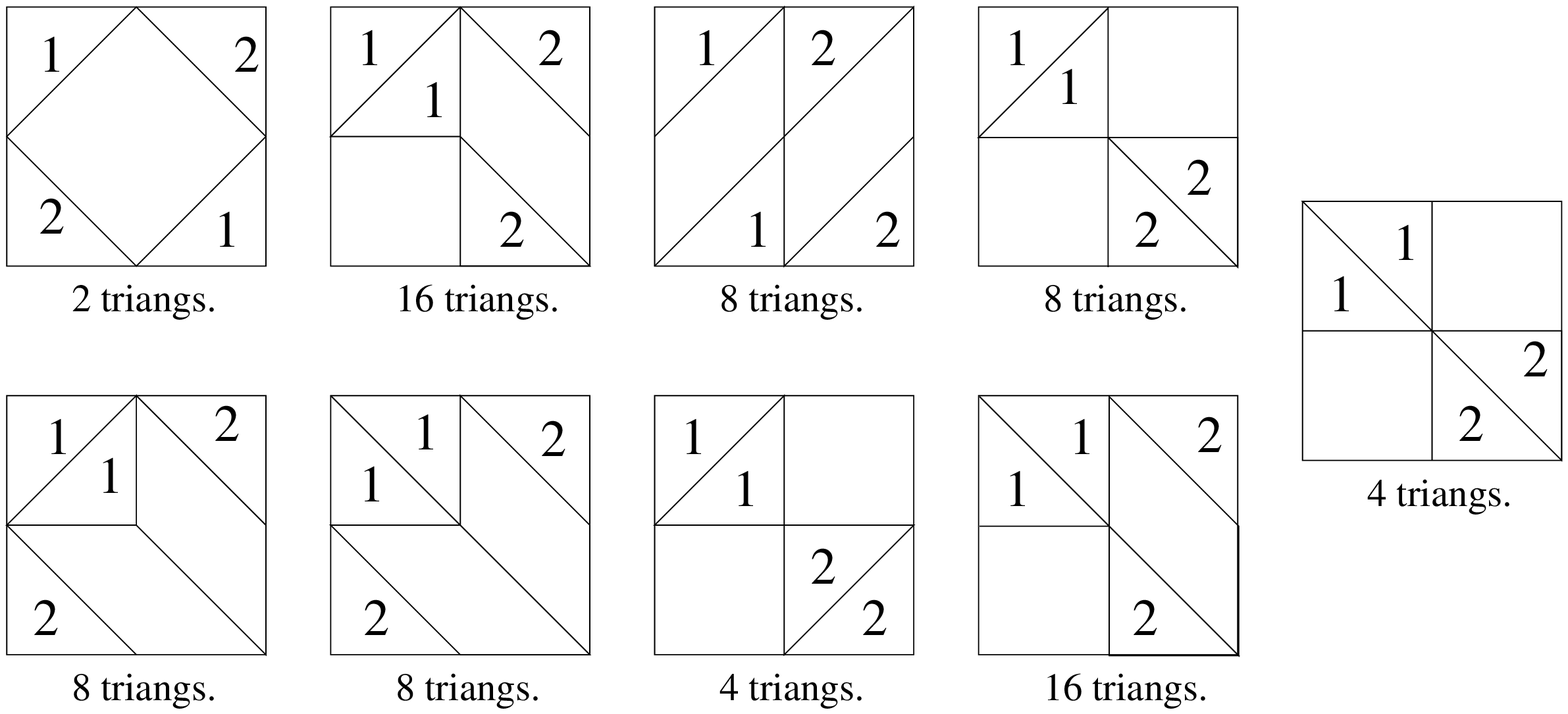}
\end{center}
\caption{The 74 triangulations of a regular cube}
\label{fig:3cube}
\end{figure}

Figure~\ref{fig:3cube} illustrates this result. The figure shows all
the triangulations of the 3-dimensional cube (the Cayley
embedding of two equal squares), pictured as mixed subdivisions.
Only one representative modulo the
symmetries of the square and modulo the exchange of the labels $1$ and
$2$ is shown. In each picture, the two triangles labeled $1$ are the
triangulation of the bottom square and the two triangles labeled $2$
are the triangulation of the top square of the cube. 
Knowing that information is
enough to recover the mixed subdivision labeling (and hence, 
the corresponding triangulation of the 3-cube).

It may be of
interest also to define {\em pure mixed subdivisions} as those for which
every Minkowski cell $B=\sum B_i$ satisfies $\dim(B)=\sum
\dim(B_i)$. If all the $P_i$'s are simplices then pure is
equivalent to fine. In general, pure
mixed subdivisions form a lower ideal in the poset of mixed subdivisions and
contain all the fine ones. Their Minkowski cells are (combinatorially)
products of their summands. Corollary~\ref{coro:labels} and
Theorem~\ref{thm:free} hold
for pure mixed subdivisions, with the same proofs.

\subsection{The product of two simplices}
The product $\Delta^p \times \Delta^q$ of two simplices can be
considered a Cayley embedding in two ways: $q+1$ copies of $\Delta^p$
or $p+1$ copies of $\Delta^q$. 
It is easy to conclude from the previous results that in this case the
``Minkowski labeling'' of a subdivision can totally be
neglected.

\begin{lemma}
\label{lemma:sumofsimplices} Let $B\subset \reals^d$ 
be a polytope and $\Delta\subset\reals^d$ be a
$d$-simplex.  If there is a way of writing $B$ as a Minkowski sum of
(positive dimensional, not necessarily distinct) faces of $\Delta$,
then this way is unique, modulo reordering.
\end{lemma}

\begin{proof} 
We argue by induction on the dimension of $\Delta$. The case where
$\Delta$ is a segment is trivial: $B$ is the sum of as many copies of
$\Delta$ as indicated by its length.

For the inductive step, let $\Delta^0$ be a facet of $\Delta$, and let
$p_0$ be the opposite vertex. For every polytope $P$ in
$\reals^d$ let $P^0$ denote the face of $P$ in the direction
of $\Delta^0$ (that is, the face containing the outer normal vector to
$\Delta^0$ in its relatively open outer normal cone)
and let $P^1$ denote the face in the opposite direction (that is, the
face whose outer normal cone contains the opposite vector).

Since $B$ can be written as a Minkowski sum of faces of $\Delta$,
every face of it can. In particular, the faces $B^0$ and $B^1$. 
Moreover, any decomposition $B=F_1+\cdots+F_k$ will restrict to decompositions 
$B^0=F^0_1+\cdots+F^0_k$ and $B^1=F^1_1+\cdots+F^1_k$.
By inductive hypothesis, the decompositions of $B^0$ and $B^1$ are
unique, so we can assume we know them, except
perhaps for the ordering and for the fact that some $F^\alpha_i$ may
be points and we cannot recover them.

 For every $i$,  $F^1_i$
must be either equal to $F^0_i$ or to the point $p_0$. This allows us to
assume that we have matched each $F^0_i$ to its $F^1_i$ and that we
can recover $F_i$, with the only
exception of the summands where both $F^0_i$ and $F^i_1$ 
are points (and hence $F_i$ is a
segment containing $p_0$). But after we subtract from $B$ all the
summands which are not segments, what remains is a parallelotope
(Minkowski sum of linearly independent segments) whose Minkowski
decomposition is straightforward, and unique.
\end{proof}

\begin{theorem}
\label{thm:sumofsimplices}
Let $\Delta$ be a simplex.
Let $S$ be a polyhedral subdivision of $k\Delta$ and assume that every
cell $B$ can be written as a Minkowski sum of faces of $\Delta$.
Then:
\begin{enumerate}
\item $S$ can be labeled as a mixed subdivision.
\item The labeling is unique, modulo reordering
of the $k$ summands.
\end{enumerate}
\end{theorem}

\begin{proof}
Uniqueness follows immediately from Theorem~\ref{thm:labels} and Lemma~\ref{lemma:sumofsimplices}. 
Indeed, the $S_i$'s are known because the
simplex has only the trivial subdivision, and which cells collapse 
to full-dimensional in $S|_{\{i\}}$ is also known (modulo reordering
of the factors): those which have a full-dimensional summand in their unique
decompositions as Minkowski sum of faces of $\Delta$.

To prove existence, let us show how to find the $i$th summand of all
cells of $S$. First, identify the cells that have the full simplex
$\Delta$ as one of the summands in their (unique, by Lemma~\ref{lemma:sumofsimplices}) Minkowski decomposition. There will be
exactly $k$ such cells, counted with multiplicity if some have
$\Delta$ as a repeated summand. Assign the numbers $1$ to $k$ to them
arbitrarily.

Once this is done, fix an index $i\in\{1,\dots,k\}$. As in the previous
results, starting from that cell we can conclude what the $i$th
summand of every other cell should be, starting with those where this
summand has codimension 1, then 2, etc. The only difficulty is to show
that this assignment is globally consistent, meaning that the $i$th
summand obtained for a given cell $B$ is independent of the path that
led from the cell with a full-dimensional $i$th summand to $B$.

To show consistency we use the 
following idea: assume without loss of generality that
$i=k$, and in each cell of $S$ shrink by a factor of $\lambda$
the $k$-th summand obtained by the previous method. Since being a polyhedral
subdivision is a local property, this produces a polyhedral
subdivision of $(k-1)\Delta +\lambda \Delta$. In the limit where
$\Delta=0$ we get a polyhedral subdivision of $(k-1)\Delta$ all of
whose cells are Minkowski sums and, by induction on $k$, the labelings
are consistent (and unique).
\end{proof}

\section{Mixed subdivisions of $k\Delta^2$ and lozenge tilings}
\label{sec:lozenge1}
The triangulations of $\Delta^1 \times \Delta^{k-1}$ are
well-understood. Their number is $k!$ and they form one only orbit
under the action of the affine symmetry group $S_k$ of $\Delta^{k-1}$
(see~\cite[Section 7.3.C]{GKZbook}). All this can be easily derived
from the Cayley Trick: The Minkowski sum of $k$ copies of a segment
is just a segment $k$ times longer, and the fine mixed subdivisions
of it are the $k!$ ways of placing the segments one after
another (in other words, the fine mixed subdivisions differ only by
the ``labeling"). This and the next section are devoted to the next case,
$\Delta^2 \times \Delta^{k-1}$.

By the Cayley Trick, subdivisions (resp., triangulations) of
$\Delta^{2} \times \Delta^{k-1}$ are in bijection with the mixed
subdivisions (resp., fine mixed subdivisions) of the Minkowski sum of
$k$ copies of the triangle $\Delta^2$. Let us denote $T_k$ this
Minkowski sum, which we think of as an equilateral triangle of size
$k$. Hence $T_1=\Delta^2$ is a triangle of unit size, whose vertices
we denote $a$, $b$ and $c$.

\subsection{Fine mixed subdivisions of $T_k$}
\begin{figure}[htb]
\begin{center}
\includegraphics[width=2.2in]{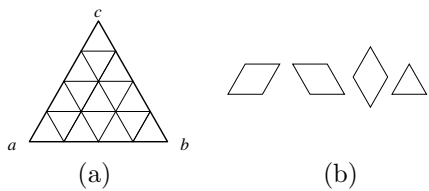}

 (a) \hskip 2.7 cm (b) 
\caption{(a) The triangular tiling of $T_4$. (b) The tiles allowed in
a lozenge tiling of $T_k$}
\label{fig:tiling}
\end{center}
\end{figure}

Fine mixed cells must be sums of faces from the summand triangles and
the sum of dimensions of the faces involved must be 2. This leaves two
possibilities: either one of the $k$ triangles plus vertices from the
other ones or the sum of two non-parallel edges of two triangles plus
vertices from the other ones. In other words, they are the upward
triangles and lozenge tiles in the next definition:

\begin{definition}
  Let $T_{k}$ be an equilateral triangle in the plane with side length $k$.
  Consider it tiled into $k^2$ equilateral triangles of side length 1,
  $(k^2+k)/2$ of them parallel to $T_k$ (we call them {\it upward triangles})
  and $(k^2-k)/2$ of them opposite to $T_k$ ({\it downward
    triangles}). See Figure~\ref{fig:tiling} where $k=4$.

 A {\it lozenge} in $T_k$ is the union of a pair of adjacent
  triangles in the tiling (one upward and one downward). A {\it lozenge
    tiling} (or {\em rhombus tiling}) 
of $T_k$ is a decomposition of $T_k$ into $(k^2-k)/2$ lozenges and
  $k$ upward triangles.
\end{definition}

The terms ``lozenge'' and ``rhombus''  tilings are synonyms in the
literature but they are used normally in a sense different from ours;
they refer to tilings by only lozenges, of a shape containing as many
upward as downward triangles.
For example, a classical object of study is the set of
lozenge tilings of the centrally symmetric hexagon with sides 
of lengths $a$,
$b$, $c$, $a$, $b$ and $c$. They are in bijection with plane
partitions that fit into an $a\times b\times c$ box. A classical
result of MacMahon gives the number of them. See, for
instance,~\cite{Propp}, 
for more information on this subject.
Theorem~\ref{thm:sumofsimplices} implies:

\begin{theorem}
\label{thm:lozenge-mixed}
Every lozenge tiling of $T_k$ 
admits a labeling as a mixed subdivision of $k$
copies of a triangle. Moreover, in any such a labeling,
\begin{enumerate}
\item Each of the $k$ copies of $\Delta^2$ appears exactly once as a
summand in one of the upward triangles of the tiling.

\item Specifying an assignment of the $k$ copies of $\Delta^2$ to the
$k$ upward triangles uniquely determines the labeling of the
mixed subdivision.
\end{enumerate}
\end{theorem}

%

\begin{figure}[htb]
\begin{center}
\includegraphics[height=2.2cm]{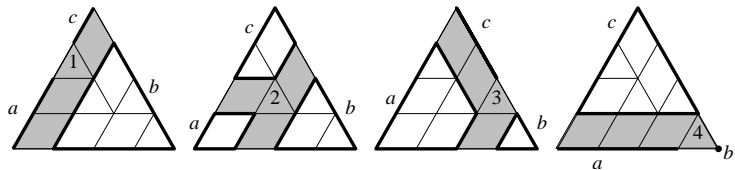}
\caption{The four zones in a certain lozenge tiling of $T_4$}
\label{fig:zones}
\end{center}
\end{figure}

There is a more direct way of proving Theorem~\ref{thm:lozenge-mixed},
that explicitly tells how to get the labelings. First, a simple
counting argument shows that there are
exactly $k$ upward triangles in every lozenge tiling: in the triangular tiling of $T_k$ there are $k$
more upward than downward triangles. 
After assigning the numbers $1$ to $k$ to
the $k$ upward triangles arbitrarily, we can define the $i$th {\em zone} of
the lozenge tiling as the union of the $i$th upward triangle plus the
lozenges that are obtained from it by parallel sweep of its
three edges along the tiling. Figure~\ref{fig:zones} shows the four
zones in a certain lozenge tiling of $T_4$. The name ``zone'' is borrowed from
a somewhat similar concept in zonotopal tilings. 

The complement of (the relative interior of) any of the zones consists of
three regions, one containing each of the three vertices of $T_k$. We
label the three regions as $a$, $b$ or $c$  depending on the
vertex of $T_k$ they contain.  Every lattice point of $T_k$ lies in exactly one
of the three regions. Then, the $i$th summand of a cell $B$ in the
tiling equals the convex hull of the vertices of $T_1$ corresponding to the
(closed) regions of the complement of the $i$th zone 
intersected by $B$.

As an example,
the central triangle in figure~\ref{fig:zones} gets labeled as
\[
\{b\} + \{a,b,c\}+ \{a\}+ \{c\}
\]

%
%
\begin{corollary}
\label{coro:lozenge-mixed}
There is a bijection between lozenge tilings of $T_k$ and orbits of
triangulations of $\Delta^2\times \Delta^{k-1}$ by the action of the
symmetry group of $\Delta^{k-1}$. In particular, the number of
triangulations of $\Delta^2\times \Delta^{k-1}$ equals $k!$ times the
number of lozenge tilings of $T_k$.
\end{corollary}

\begin{figure}[htb]
\begin{center}
\includegraphics[width=3.8in]{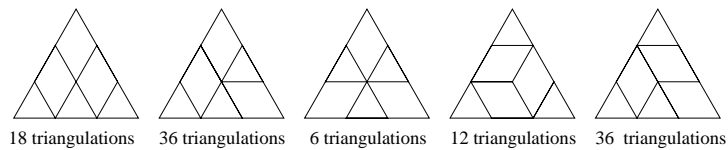}
\caption{The 108 triangulations of $\Delta^2\times \Delta^2$}
\label{fig:twobytwo}
\end{center}
\end{figure}

As a first application of this result, we can ``draw'' all the
triangulations of $\Delta^2\times\Delta^2$, that is to say, all mixed
subdivisions of $T_3$. Modulo the symmetries
of the triangle, there are 5 lozenge tilings of $T_3$, 
displayed in Figure~\ref{fig:twobytwo}. Each represents as many $S_3$-orbits of
triangulations as lozenge tilings in its orbit modulo the 
symmetries of the triangle\footnote{Observe that in this sentence there are
of two different $S_3$-actions: one on the labels, the other on $T_3$.
They correspond to the symmetries in the two factors of
$\Delta^2\times \Delta^2$.}.
Hence, the number of lozenge tilings of $T_3$ is
$3+6+1+2+6=18$ and the number of triangulations of
$\Delta^2\times\Delta^2$ is $18\cdot 6 = 108$.  The reader can compare
Figure~\ref{fig:twobytwo} with Figure 39 in~\cite{GKZbook} (page 150),
where a different representation of the triangulations of
$\Delta^2\times \Delta^2$ is used. There, the vertices of
$\Delta^{k-1}\times \Delta^{l-1}$ are represented as a $k\times l$
grid, and each simplex of a triangulation is represented by marking
some squares in the grid.  Incidentally, comparing the two figures
the reader can easily detect an error in the adjacency graph of one of
the triangulations of~\cite[Fig. 39]{GKZbook}.

\subsection{Non-fine mixed subdivisions of $T_k$}
What should a non-fine mixed subdivision look like?  In the first
place, it must be a polyhedral subdivision in the usual sense, with
vertex set contained in the ${k+2 \choose 2}$ lattice points of
$T_k$. Second, it must possess fine mixed refinements, hence each cell
must be a convex union of lozenge tiles and upward triangles. We
consider any such union as a hexagon with three pairs of parallel
edges, each pair parallel to one edge of $T_k$. 
Hexagons may degenerate to have some edges of length zero.

\begin{lemma}
\label{lem:mixed-nonfine}
Let $B$ be a convex union of triangles of the triangular tiling. Then,
the following properties are equivalent:
\begin{enumerate}
\item $B$ is a Minkowski sum of faces (perhaps with repetition) of the
unit upward triangle.
\item $B$ can be tiled by lozenges and upward triangles.
\item $B$ contains at least as many upward triangles as downward triangles.
\item In each of the three pairs of parallel sides of $B$ the one with
the same outer normal as a side of the unit upward triangle is at
least as long as the opposite one.
\end{enumerate}
\end{lemma}

\begin{proof}
If $B$ has a pair of opposite sides of positive length, then reducing
both by one unit does not affect whether $B$ satisfies any of the
properties. Hence, there is no loss of generality in assuming that $B$
has no pair of opposite sides, hence it is a triangle (or a
point). If it is an upward triangle of any size (or a point) then the four
conditions hold; if it is a downward triangle then none of them does.
\end{proof}

As in the fine case, Theorem~\ref{thm:sumofsimplices}
has the following consequence. Observe, however, that the orbits of the
permutation group may now have cardinality smaller than $k!$ (an extreme
example of this is the trivial subdivision).

\begin{theorem}
\label{thm:mixed-nonfine}
A polygonal subdivision of $T_k$ can be labeled as a mixed subdivision
if and only if all the cells are in the conditions of Lemma~\ref{lem:mixed-nonfine}.
The labeling is unique modulo the action of the
permutation group.
\end{theorem}

A direct way of getting the Minkowski labeling in this case is as follows:
We define the excess of a cell as the difference between upward
and downward triangles contained in it. Our assumption is that all
cells have non-negative excess. Actually, the excess of a cell has an
interpretation in any of the four settings of Lemma~\ref{lem:mixed-nonfine}: it is the number of upward triangles in a
Minkowski decomposition, the number of upward triangles in a lozenge
tiling, and the difference in length between any of the three pairs of
opposite edges.

The total excess in the tiling is clearly $k$. Let
us distribute the numbers 1 through $k$ to the different cells, giving
a cell as many numbers as its excess. In much the same way as we did
for lozenge tilings, we can define the $k$ zones of the polyhedral
subdivision: The $i$-th
zone contains the cell to which we assigned the label $i$ 
and then three arms, obtained as the cells adjacent to
it in the directions towards the three edges of $T_k$, and then the
ones adjacent to these, and so on. The main difference with the zones
in a lozenge tiling is that now, as we travel along an arm, the edges
that we cross may increase in length from one cell to the next. Also,
the $i$th and $j$th zones may coincide (if the $i$th and $j$th excess
reside in the same cell) or an arm of one zone be contained in the
other zone.

As in the case of lozenge tilings, the definition of $i$-th zone
classifies cells into seven types: the one labeled $i$, the ones in the
three arms and the ones in the three regions of the complement of the
zone (some of the last six types may be empty). This classification
says whether the $i$-th summand in the mixed cell expression of a
given cell is going to be $\{a,b,c\}$ (if the cell is labeled $i$),
$\{a,b\}$, $\{a,c\}$ or $\{c,b\}$ (if the cell lies in one of the three arms, depending on the
edge of $T_k$ they are heading to) or just $\{a\}$, $\{b\}$, or
$\{c\}$ (if the cells in the complement of the zone).

As an example, the right part of Figure~\ref{fig:nonregular} shows a
valid polygonal subdivision of $T_8$. The cells labeled 1 to 8 have
excess 1. The two unlabeled cells have excess zero.  The dots, arrow
and shading in the figure are there for later use. In order
to illustrate the above concepts, Table~\ref{table:tencells}
shows the mixed subdivision labeling (i.e., the Minkowski decomposition)
of the ten cells. Cells 1 to 8 appear first in the
list, then the shaded hexagon and finally the parallelogram. The
eight columns of summands correspond to the eight zones.
\begin{table}
\[
\begin{array}{c}
\{a,b,c\} + \{a,c\}\ +\ \{a,c\}\ +\ \ \{a\}\ \ +\ \ \{a\}\ \ + \ \ \{a\}\ \  +\ \{a,b\}\ +\ \ \{a\}. \medskip \\
\ \ \{b\}\ \ + \{a,b,c\} + \ \{a,c\}\ +\ \ \{a\}\ \ +\ \ \{a\}\ \ +\ \ \{a\}\ \ +\ \ \{b\}\ \ +\  \{a,b\}.\medskip \\
\ \ \{b\}\ \ +\ \ \{b\}\ \ + \{a,b,c\} +\ \ \{a\}\ \ +\ \{a,b\}\ +\ \ \{a\}\ \ +\ \ \{b\}\ \ +\ \ \{b\}.\medskip \\
\ \ \{b\}\ \ + \ \ \{b\}\ \ +\ \{b,c\} \ + \{a,b,c\} +\ \ \{b\}\ \   +\  \{a,b\} \ + \ \ \{b\}\ \   + \ \ \{b\}.\medskip \\
\ \ \{b\}\ \ +\  \{b,c\}\  +\ \ \{c\}\ \ + \ \{a,c\}\  + \{a,b,c\} + \ \ \{a\}\ \  + \ \ \{b\}\ \  + \ \ \{b\}.\medskip \\
\ \{b,c\}\  +\ \ \{c\}\ \  + \ \ \{c\}\ \  + \ \ \{c\}\ \  + \ \ \{c\}\ \  + \{a,b,c\} + \ \ \{b\}\ \ +\ \ \{b\}.\medskip \\
\ \ \{c\}\ \  + \ \ \{c\}\ \  + \ \ \{c\}\ \  +\  \{a,c\}\  + \ \{a,c\}\  + \ \ \{a\}\ \   + \{a,b,c\} + \ \ \{a\}.\medskip \\
\ \ \{c\}\ \  + \ \ \{c\}\ \  + \ \ \{c\}\ \  + \ \ \{c\}\ \  + \ \ \{c\}\ \  + \ \{a,c\}\  + \ \{b,c\}\  + \{a,b,c\}.\medskip \\
\ \{bc\}\  + \ \ \{c\}\ \  + \ \ \{c\}\ \  +\ \{a,c\}\ +\ \{a,c\} \ + \ \ \{a\}\ \  + \ \ \{b\}\ \ + \ \{a,b\}. \medskip \\
\ \ \{b\}\ \  + \ \{b,c\}\  + \ \ \{c\}\ \  + \ \ \{c\}\ \  + \ \{b,c\}\  + \ \{a,b\}\  + \ \ \{b\}\ \ + \ \ \{b\}.\\
\end{array}
\]
\medskip
\caption{Minkowski decomposition of the cells in the 
non-regular subdivision of $\Delta^2\times \Delta^7$ of
Figure~\ref{fig:nonregular} (right)}
\label{table:tencells}
\end{table}

\section{Subdivisions of $\Delta^2\times \Delta^{k-1}$}
\label{sec:lozenge2}

In this section we list several properties of triangulations and
subdivisions of  $\Delta^2\times \Delta^{k-1}$ which can be derived
from representing them as lozenge tilings.

\subsection{Non-regular subdivisions of $\Delta^2\times \Delta^{k-1}$}

In order for a mixed subdivision to be coherent it has first to be
regular as a subdivision in the standard sense (i.e., the projection
of the lower hull of a polytope in one dimension more).
Hence:

\begin{figure}[htb]
\begin{center}
\includegraphics[height=2in]{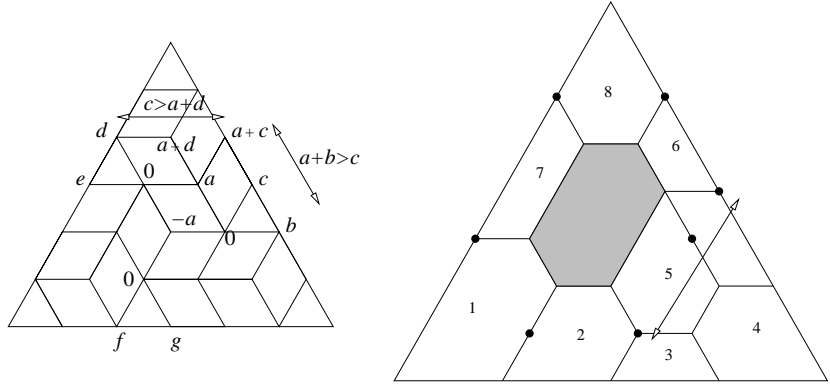}
\end{center}
\caption{A non-regular triangulation of $\Delta^5\times \Delta^2$ (left)
and a coarse non-regular subdivision of $\Delta^2\times \Delta^7$
(right)}
\label{fig:nonregular}
\end{figure}

%
%

\begin{proposition}
\label{prop:nonregular}
The tilings of Figure~\ref{fig:nonregular} represent 
a non-regular triangulation of $\Delta^5\times \Delta^2$
and a coarse non-regular subdivision of $\Delta^2\times \Delta^7$,
respectively.
\end{proposition}

\begin{proof}
The proof of non-regularity is sketched in the picture in both
 cases. In the left, if a lifting existed, there would be no loss
 of generality (by addition of an affine function to all the heights)
 in assuming that the three neighbors of the central point get height
 zero. The central point must then get a negative height that we
 denote $-a$.  We let $b,c,d,e,f$ and $g$ denote the heights of
 certain boundary points, as shown in the figure. From these heights
 some others can be deduced, and in particular the figure shows how to
 conclude that $b>d$. The same arguments applied cyclically show that
 $d>f$ and $f>b$, which is impossible.  

For the picture on the right, there is no loss of generality in
 assuming height zero for all the vertices of the shaded cell. Then,
 the seven marked points can easily be proved to get all the same
 height, but this contradicts convexity at the edge between region 5
 and its adjacent parallelogram.  That the right picture represents a
 coarse subdivision follows from the fact that it is coarse as a
 subdivision of $T_8$ in the standard sense.
\end{proof}

Non-regular triangulations of the product of two simplices were first
constructed by de Loera~\cite{Loera}, for $\Delta^3\times\Delta^3$. 
He also proved that all triangulations of
$\Delta^2\times\Delta^k$ are regular, up to $k=4$. Later,
Sturmfels~\cite{Sturmfels-book}, constructed a non-regular
triangulation of 
 $\Delta^2\times\Delta^5$, hence concluding that 
$\Delta^k\times\Delta^l$  has non-regular triangulations if and only
 if $(k-1)(l-1)\ge 4$. In particular, the non-regular lozenge tiling of $T_6$
 that we show is smallest possible. As for our second example, to the
 best of our knowledge it is the first known coarse non-regular
 subdivision of the product of two simplices.  Observe that coarse
 subdivisions of polytopes in general, and of products of simplices in
 particular, are not well-understood objects.

\begin{figure}[htb]
\begin{center}
\includegraphics[height=2.5in]{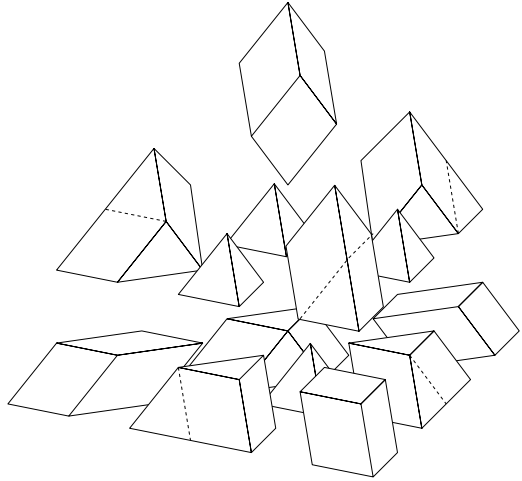}
\end{center}
\caption{A regular subdivision of $\Delta^3\times \Delta^3$ and a 
non-regular refinement of it}
\label{fig:3x3nonregular}
\end{figure}

The Cayley Trick can also be used to picture non-regular
triangulations of $\Delta^3\times \Delta^3$.
Figure~\ref{fig:3x3nonregular} is our attempt to do so.
The picture shows (an explosion of)
the coherent mixed subdivision of $4\Delta^3$
produced by the following lifting matrix. Each row represents the
lifting of one of the four copies of $\Delta^3$:
\[
\left(
\begin{matrix}
0 & 1 & 1 & 1 \\
1 & 0 & 1 & 1 \\
1 & 1 & 0 & 1 \\
1 & 1 & 1 & 0 \\
\end{matrix}
\right).
\]
This subdivides $4\Delta^3$ into 14 cells: four parallelepipeds in the four
corners of $4\Delta^3$; four tetrahedra incident to 
the center of the four facets
of $4\Delta^3$; and six ``Minkowski sums of two triangles'' along the six
edges of $4\Delta^3$. Another way of describing this subdivision is
that it is obtained by cutting $4\Delta^3$ with the four planes
through its centroid and parallel to its facets.

\begin{figure}[htb]
\begin{center}
\includegraphics[height=0.95in]{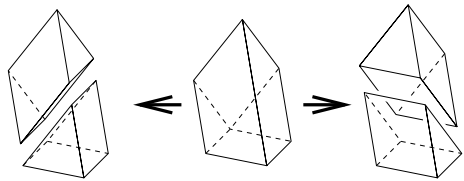}
\end{center}
\caption{Refining the non-fine cells of Figure~\ref{fig:3x3nonregular}}
\label{fig:3x3flip}
\end{figure}

This mixed subdivision is not fine, because the six special cells
along edges of $4\Delta^3$ can be refined into two triangular prisms
each, as Figure~\ref{fig:3x3flip} shows. If the six
cells are refined in the particular ``skew'' way sketched by dashed
lines in Figure~\ref{fig:3x3nonregular}, then the fine mixed
subdivision obtained is not coherent. (It is actually a non-regular
subdivision; the proof is easy and left to
the reader).  Hence, it corresponds to a non-regular triangulation of
$\Delta^3\times \Delta^3$.

It is interesting to observe that $\Delta^{k-1}\times\Delta^{l-1}$ has non-regular
triangulations if and only if there is a matroid on $k+l$ elements and of rank $k$
which is not representable over the reals. This was first noticed in~\cite{Sturmfels-book},
where the subdivisions in Figure~\ref{fig:3x3nonregular} and the left part of
Figure~\ref{fig:nonregular} were related, respectively, to the Vamos and the non-Pappus 
matroids.

\subsection{Lozenge flips versus bistellar flips}
A basic concept to understand the set of all triangulations of a
polytope is that of geometric bistellar flip. Roughly speaking, it is
the minimum possible difference (the ``elementary move'') between two
triangulations. One simple definition, (see ~\cite{Santos-refine}) is that two triangulations differ by a bistellar
flip if and only if they are the only two refinements of a certain
polyhedral subdivision. A more explicit definition that says what the
difference between the two triangulations has to be for this to happen
is contained, for
example, in~\cite{GKZbook,Triangbook}. We do not need it here.

\begin{figure}[htb]
\begin{center}
\includegraphics[height=1in]{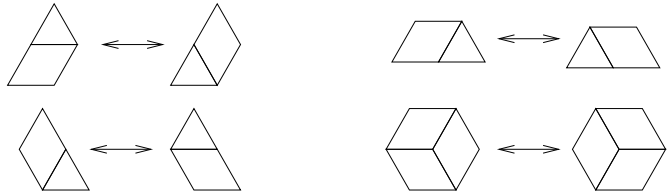}
\end{center}
\caption{The four types of lozenge flips}
\label{fig:flips}
\end{figure}

\begin{definition}
\label{defi:lozenge-flip}
We say that two lozenge tilings of $T_k$ differ by a {\em lozenge
flip} if one can be obtained from the other by one of the four
substitutions of tiles shown in Figure~\ref{fig:flips}.  More
precisely, the first three will be called {\em trapezoid flips}
and the last one a {\em hexagon flip}.
\end{definition}

\begin{proposition}
\label{prop:flips}
Let $L$ and $L'$ be two labeled lozenge tilings of $T_k$,
corresponding to two triangulations $S$ and $S'$ of $\Delta^2\times
\Delta^{k-1}$. Then, $S$ and $S'$ differ by a bistellar flip if and
only if the following three properties hold:
\begin{enumerate}
\item $L$ and $L'$ differ by a lozenge flip.
\item The labeling of upward triangles is the same in $L$ and $L'$
(except for the displacement of the triangle in case of a trapezoid
flip).
\item For a trapezoid flip, the arm 
of the triangle affected by the flip in the direction of the big edge
of the trapezoid does not change
with the flip (it is only translated).
\end{enumerate}
\end{proposition}

\begin{proof}
Let us look at what a polyhedral subdivision of $T_k$ has to look like
in order to admit only two lozenge refinements. First, all cells must
be convex and individually admit only two lozenge refinements. The
possibilities are a hexagon as the one in a hexagon flip, a trapezoid,
or the union of two parallel lozenges with a common edge. If a hexagon
arises, then its refinements are independent of the refinement of any
other cell, which means that no other refinable cell can be
present. Hence, the two tilings differ by a hexagon-flip. If a
trapezoid or union of two lozenges arises, then 
the edge (or edges) of length two in that cell must
be propagated up to the boundary of $T_k$ on one side and to a
trapezoid on the other side. The flip is a trapezoid flip and
satisfies the arm condition in the statement. That the labels must be
the same in the two lozenge tilings is trivial.
\end{proof}

\begin{figure}[htb]
\begin{center}
\includegraphics[height=2.7in]{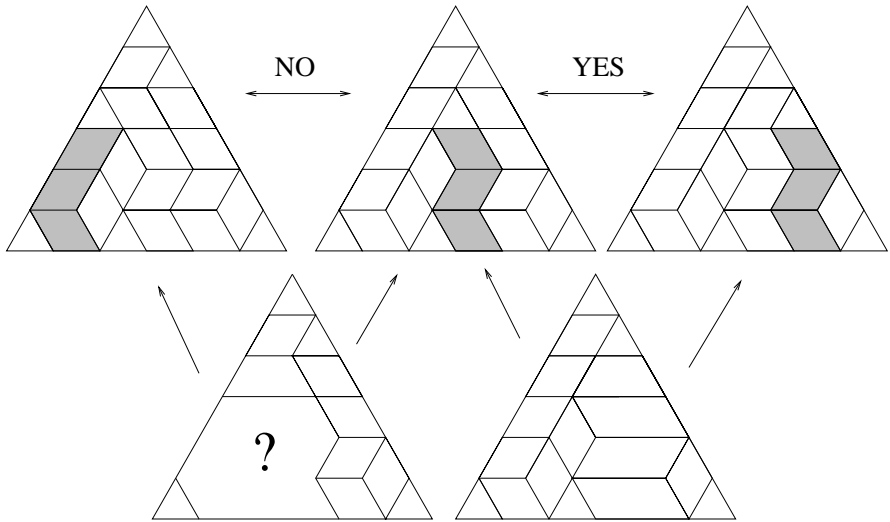}
\end{center}
\caption{A trapezoid flip is not always a bistellar flip}
\label{fig:noflip}
\end{figure}

As an example, the two lozenge tilings on the right part of Figure~\ref{fig:noflip} represent triangulations which differ by a bistellar
flip, because they are the two refinements of the subdivision in the
bottom-right. The two lozenge tilings on the left differ by a lozenge
flip, but not by a bistellar flip.

\begin{theorem}
\label{thm:connected}
\begin{enumerate}
\item The set of all (labeled) lozenge tilings of $T_k$ is connected
under trapezoid flips.
\item The set of triangulations of $\Delta^2\times\Delta^{k-1}$
is connected under geometric bistellar flips.
\end{enumerate}
\end{theorem}

\begin{figure}[htb]
\begin{center}
\includegraphics[height=1.3in]{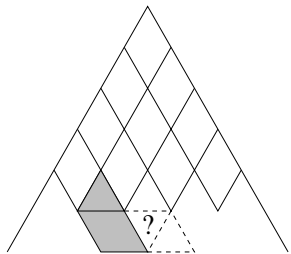}
\end{center}
\caption{Every lozenge tiling has {\em some}
 trapezoid flip which is a bistellar flip}
\label{fig:goodflip}
\end{figure}

\begin{proof}
Any triangle not on the bottom row of a lozenge tiling is adjacent to
a lozenge below it.  Performing the trapezoid flip there produces a
triangle one level lower. With this idea, we can eventually arrive to
a lozenge tiling with all its triangles on the bottom row, and there
is only one such tiling. We have not taken care of labels, but once we
have the tiling with all triangles in the bottom row there is a
sequence of three trapezoid flips which exchanges two consecutive
triangles. Hence, any permutation of the labels can be implemented as
a sequence of trapezoid flips, too. This proves part (1).

For part (2), we proceed similarly. Not all trapezoid flips are
bistellar flips, but we can prove that unless all the triangles are in
the bottom row there must be some trapezoid flip which decreases the
height of the triangle involved and which is a geometric bistellar
flip. To see this, start a triangle $i$ with maximum height in the tiling.
This implies that all the rows above that triangle are tiled with
vertical lozenges, as in Figure~\ref{fig:goodflip}. Let us consider
the trapezoid flip that would decrease the height of $i$ and,
more specifically, at the downward looking triangle next to the
lozenge below $i$, in the direction of the big side of the trapezoid.
There are two possibilities for the lozenge containing that downward
triangle: if it is a vertical lozenge, then the trapezoid flip at $i$
is a geometric bistellar flip, and we are done. If it is not, then 
above it there is another triangle at the same height as $j$, and we
get a trapezoid flip looking in the same direction and ``closer to the
boundary''. Repeating this process we must eventually arrive at a
trapezoid flip which is itself a bistellar flip.
\end{proof}

Further analysis of the above proof gives bounds on the
diameter:

\begin{corollary}
\label{coro:diameter}
The graph of lozenge tilings of $T_{k+1}$ and the graph of triangulations
of $\Delta^2\times\Delta^k$ both have diameter in $\Theta(k^2)$.

More precisely, the graph of unlabeled lozenge tilings has diameter at
least ${k\choose 2}$ and the graph of triangulations of
$\Delta^2\times\Delta^k$ has diameter at most  $5{k\choose 2}$.
\end{corollary}

\begin{proof}
For a quadratic lower bound, observe that at least ${k\choose 2}$
lozenge flips are needed to go from the lozenge tiling with all
triangles on the bottom to the lozenge tiling with all triangles on
one side. This is so because each lozenge flip changes the height of
only one triangle and only by one unit. 

For an upper bound, the process in the proof of Theorem~\ref{thm:connected}
shows how to go from any lozenge tiling to the one with all triangles
in the bottom by a sequence of $\sum(\height_T(i))\le {k\choose 2}$ 
bistellar flips. As
mentioned there, once we are in that lozenge tiling we can permute any
two labelings with three times the number of pairs of indices which
are ordered differently, that is, at most $3{k \choose 2}$ bistellar
flips. With another  ${k \choose 2}$ flips we can go back to the
second lozenge tiling.
\end{proof}

The constants in the previous statement can surely be improved. For
example, instead of going to the tiling with all triangles on
the bottom, we can choose to go to the
tiling with all triangles on one side, which provides a
different definition of height. For each triangle, the sum of its
three heights is clearly $k-1$, so that with respect to one of the three
sides we get 
\[
\sum(\height_T(i)) + \sum(\height_{T'}(i)) \le 2 k(k-1)/3
\]
instead of the $2{k \choose 2}=k(k-1)$ used in the proof.


\begin{question}
Is there a triangulation of $\Delta^{2}\times\Delta^{k-1}$
with less than $2k-2$ bistellar flips? Observe that $2k-2$ is the
dimension of the corresponding secondary polytope, hence it is  a
lower bound for the number of flips of every regular
triangulation. Also, it is easy to prove that every lozenge tiling has
at least $2k-2$ lozenge flips: The $k$ upward
triangles have a total of $3k$ sides and at most $k+3$ of them are in
the boundary. This implies there are at least $2k-3$ trapezoid
flips and that there are exactly that number if and only if three
triangles are at the corners and the others are on the boundary. In
this case, the lozenges produce a lozenge tiling of a simply
connected region in the standard sense. Since lozenge tilings of a
simply connected region are connected by hexagon flips, there has to
be at least one hexagon flip.
\end{question}

\subsection{Counting lozenge tilings}

The number of lozenge tilings of $T_{k}$ can be computed in the
following way. Let $k$ be fixed, and let $S$ be a subset of
$\{1,2\dots,k\}$. We classify the lozenge tilings of $T_k$ according
to what triangles they have in the bottom line. More precisely, let
\begin{itemize}
\item $f_k(S)$ denote the number of lozenge tilings of $T_k$ which have triangles exactly in the positions of the bottom line given by $S$.
\item $g_k(S)$ denote the number of lozenge tilings of $T_k$ which have triangles at least in the positions of the bottom line given by $S$.
\end{itemize}
Clearly,
\begin{equation}
g_k(S) = \sum_{S\subseteq S'} f_k(S').
\label{eq:gS}
\end{equation}
But, moreover,
\begin{proposition}
Let $S=\{s_1,\dots,s_j\}\ne \emptyset$, where $1\le s_1<\cdots<s_j\le k$.
If $j=1$, then $f_k(S)=g_{k-1}(\emptyset)$. If $j>1$, then:
\begin{equation}
f_k(S)=\sum_{
\begin{array}{c}
s_1\le s'_1<s_{2} \cr
\vdots \cr
s_{j-1}\le s'_{j-1}<s_{j} \cr
\end{array}
} 
g_{k-1}(\{s'_1,\dots,s'_{j-1}\})
\label{eq:fS}
\end{equation}
\end{proposition}

\begin{proof}
Between every two triangles of the bottom row there must be one and only one vertical lozenge.
Once we fixed the positions $s'_1,\dots,s'_{j-1}$ of these vertical lozenges, the ways to
complete the lozenge tiling are exactly the same as the lozenge tilings of $T_{k-1}$ containing
triangles in (at least) the positions $s'_1,\dots,s'_{j-1}$ of the bottom row.
\end{proof}

Table~\ref{table:f-g} shows all the values of $f_k(S)$ and $g_k(S)$ with
$k=1,2,3$, as well as the values of $f_4(S)$,
computed using the recursive equations~(\ref{eq:gS}) and~(\ref{eq:fS}). 
Adding all the entries of $f_4(S)$ we get the number of
lozenge tilings of $T_4$, which is $g_4(\emptyset)=187$. Hence, the
number of triangulations of $\Delta^2\times\Delta^3$ is $187\times
4!=4488$. 

\footnotesize
\begin{table}[htb]
\[
\begin{array}{|c|c|c|c|c|c|c|c|c|}
\hline
S &    \emptyset & 1 & 2 & 1,2 & 3 & {1,3} & 2,3 & 1,2,3  \cr
\hline
\hline
f_1(S) &  0           &   1    &         &            &         &        &        &        \cr        
\hline
g_1(S) &  {\bf 1}       &   1    &         &            &         &        &        &        \cr        
\hline
\hline
f_2(S) &    0   &   1    &    1    &   1    &         &        &        &        \cr
\hline
g_2(S) &   {\bf 3}   &  2    &  2    &    1    &         &        &        &        \cr
\hline
\hline
f_3(S)  & 0  & 3 & 3 & 2 & 3 & {\begin{array}{c} 2+2 \\ =4 \end{array} } 
& 2 &1\cr
\hline
g_3(S)  & {\bf 18}  & 10 & 8 & 3 & 10 & 5 & 3 & 1\cr
\hline
\hline
f_4(S) & 
0  & 18 & 18 & 10 & 
18 & {\begin{array}{c} 10+8 \\ =18 \end{array} } & 8 & 3\cr
\hline
f_4(S\cup\{4\}) & 
18  & 
{\begin{array}{c} 10+8+10 \\ =28 \end{array} } & 
{\begin{array}{c} 8+10 \\ =18 \end{array} } & 
{\begin{array}{c} 3+5 \\ =8 \end{array} } & 
10 & 
{\begin{array}{c} 5+3 \\ =8 \end{array} } & 3 & 1\cr
\hline
\end{array}
\]
\caption{The number of triangulations of $\Delta^2\times\Delta^3$ computed by hand. It
 is $4!$ times the sum of entries in the last two rows}
\label{table:f-g}
\end{table}
\normalsize

The numbers shown in Table~\ref{table:numbers} in the introduction
were computed with an 
implementation of these recursive formulas in Maple.
The computation is clearly exponential in time,
since we need to compute $2^k$ values of $f_k(S)$ and $g_k(S)$ 
for each $k$. In practice, 
the computation of each value took about five times the
previous one: 21 seconds for $k=10$ and 70 hours for $k={16}$.
By Corollary~\ref{coro:lozenge-mixed}, multiplying  the $k$th number by
$k!$ we get the number of triangulations of $\Delta^{k-1}\times
\Delta^2$. A direct approach allowed Jes\'us de
Loera and J\"org Rambau~\cite{Loera,Rambau-TOPCOM} to compute these numbers
of triangulations 
only up to $k=4$ and $k=6$ respectively.

\subsection{The asymptotic number of lozenge tilings}
\begin{figure}[htb]
\begin{center}
\includegraphics[height=1.7in]{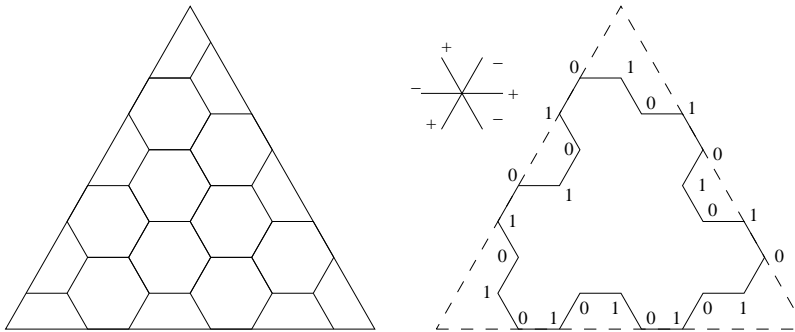}
\caption{
Proof of a quadratic lower bound for $\log(l_k)$ (left)
and a lozenge tileable region of nearly constant boundary height (right)}
\label{fig:hexagons}
\end{center}
\end{figure}

Let $l_k$ denote the number of 
lozenge tilings of  $T_k$. It is easy to show that
$l_k$ is in $e^{\Theta(k^2)}$:

\begin{itemize}
\item Since a lozenge tiling can be specified by which of the three
upward neighbors of each of the $(k^2-k)/2$ downward triangles forms a
lozenge with it, $ l_k \le 3^{(k^2-k)/2} <3^{\frac{k^2}{2}}.$

\item Assume $k$ is a multiple of 3.
$T_k$ can be tiled into $3{k/3\choose 2}=\frac{k^2-3k}{6}$ 
hexagons plus $k$ boundary
trapezoids (see Figure~\ref{fig:hexagons}), each of which can independently
be refined in two ways. Hence, $l_k\ge 2^{(k^2+3k)/6}$.

\end{itemize}

In particular, it is asymptotically not 
relevant to distinguish between labeled
and unlabeled lozenge tilings. We can as well think of
$l_k$ as the number of triangulations of $\Delta^2\times\Delta^{k-1}$.

The property that the logarithm of the number of tilings is
proportional to the area for dilations of a given shape
is well-known in the context of usual lozenge
tilings, as follows from the following result of
Cohn, Kenyon and Propp~\cite{CoKePr} (they consider mostly the case of domino
tilings, i.e., perfect matchings in a sub-region of the square grid, 
but the case of lozenge tilings arises as a particular case
in which certain edges are forbidden in the matching).
Let $R^*$ be a simply-connected region in the
plane, and let $(R_n)_{n\in \naturals}$ be a sequence of
lozenge-tileable (in the standard sense)
simply-connected regions such that 
$R_n/n$ converges to $R^*$. Let $f_n:\partial R_n\to \reals$ be the 
boundary height function of the region $R_n$ (defined below). Assume
that, after scaling it down,  $f_n/n$ converges 
to a certain function 
$f^*:\partial R^*\to \reals$. 

\begin{lemma}[Cohn et al.~\cite{CoKePr}]
\label{lemma:CoKePr}
In the above conditions, 
\begin{enumerate}
\item The logarithm of
the number of lozenge tilings of $R_n$ divided by the area of $R_n$
(measured in lozenge tiles) 
converges to a constant that depends only on $R^*$ and $f^*$.

\item This constant is maximized if $f^*=0$.
In this case it equals
\[
\frac{3}{\pi} L\left(\frac{\pi}{3}\right) \simeq 0.32306594.
\]
Here, $L(x)$ is the Lobachevsky function, defined as
\[
L(x)=-\int_{0}^{x} \log |2\sin t|\ \text{\rm d}t .
\]
\end{enumerate}
\end{lemma}

The constant in part (1) (and its specific instance in part (2))
is computed as an integral of the average (in a well-defined sense) extension 
of the boundary height function  $f^*$ to  the interior of $R^*$.

The boundary height function of a simply connected union $R$ of triangles
of the regular triangular tiling is defined as follows: choose
alternating signs for the six directions of edges in the tiling. 
Starting at
any particular boundary vertex, give height zero to that vertex and
then propagate the height along the boundary cycle of $R$,
increasing or decreasing the height by one depending on the direction
of the edge traversed. A simply connected 
 region is lozenge-tileable if and only if the 
height becomes 0 again when you return to the starting point~\cite{Thurston}.
The right part of Figure~\ref{fig:hexagons} is an 
example of a tileable region with nearly
constant boundary function.

The following statement says that 
the asymptotic entropy per unit tile is the same in our lozenge 
tilings of $T_k$ as in classical lozenge tilings of a simply connected
region with nearly constant boundary height function. The proof we give
is essentially glued from personal communications to the author by
J. Propp, H. Cohn and, specially, David Wilson:

\begin{theorem}
\label{thm:asymptotics}
\[
l_k = e^{\beta\frac{k^2}{2}\pm o(k^2)},
\]
where $\beta = 
\frac{3}{\pi} L\left(\frac{\pi}{3}\right) \simeq 0.32306594$
is the asymptotic entropy per unit tile of regions with nearly
constant boundary height function, as given by Lemma~\ref{lemma:CoKePr}.
\end{theorem}


\begin{proof}
For the lower bound, apply part (2) of 
Lemma~\ref{lemma:CoKePr} to the region
on the right of Figure~\ref{fig:hexagons}.

For the upper bound, let $f(k)$ be any  function such that,
asymptotically, $1<<f(k)<<\sqrt{k}$ (for example, 
$f(k)=k^{1/4}$). Let $S_k$ be a tiling of 
$T_k$ into a triangular grid of
about $k^2/f(k)^2$ triangles of size about $f(k)$.

For each lozenge tiling of
$T_k$, we cut $T_k$ without breaking tiles but otherwise as close as
possible to the tiling $S_k$. If a tile overlaps two cells of $S_k$,
we choose, for instance, to give that tile to the bottom of the two.
The total perimeter of the cells in $S_k$ is clearly in 
$\Theta(k^2/f(k))\subset o(k^2)$. Hence, the number of possible
ways of cutting $T_k$ produced in this way
is in $2^{o(k^2)}$.  and will not affect the final
asymptotics. Our task is to bound the number of lozenge tilings
compatible with a specific cutting. 

For this, we consider independently the cells of $S_k$ that only contain
lozenges and those that contain at least a triangle in the lozenge
tiling. Although this is not relevant for the asymptotics,
observe that which cells contain triangles (and how many of them) is
fixed by the cutting: for a specific cell, the number of triangles is
the difference between upward and downward triangles of the
triangular unit grid contained in it.

Since at most $k$ cells contain triangles, we do not need to
care much about their number of tilings. The easy argument that each
downward triangle must be matched to one of at least three upward
triangles shows that the number of tilings of each cell is 
at most $3^{f(k)^2/2}$. Hence, the cells
that contain triangles produce a factor of at most $3^{kf(k)^2/2}\in
e^{o(k^2)}$ in the final number and can be neglected.

For the cells that are tiled only with lozenges, 
we are in the situation of
Lemma~\ref{lemma:CoKePr}: the number of lozenge tilings of each
 is at most $e^{\beta \frac{f(k)^2}{2}
+o(f(k)^2)}$ and the combined number is at most
\[
\left(e^{\beta \frac{f(k)^2}{2}
+o(f(k)^2)}\right)^{k^2/f(k)^2}
=
e^{\beta \frac{ k^2}{2} +o(k^2)}.
\]
\end{proof}

\section{Tropical polytopes}
\label{sec:tropical}

Develin and Sturmfels have recently started developing the theory of
polytopes in tropical space~\cite{DevStu}. We here give a brief
account of their main results, specially in their relations to
subdivisions of the product of two simplices.

The tropical projective space of dimension $l-1$, denoted $\TP^{l-1}$,
 is the quotient of 
$\reals^{l}$ by the equivalence relation $v \sim
v+(\lambda,\dots,\lambda)$, for every $v\in \reals^l$ and every
$\lambda\in \reals$. By normalizing one of the coordinates (say the
first one) to be equal to zero we can identify $\TP^{l-1}$ to
 $\reals^{l-1}$.

The tropical hyperplane 
defined by a vector $(a_1,\dots,a_l)\in \reals^l$ is the
set of points $v\in\TP^{l-1}$ such that the minimum of the numbers
$v_i + a_i$ is achieved twice. Clearly, the hyperplanes defined
by $(a_1,\dots,a_l)$ and by $(a_1+\lambda,\dots,a_l+\lambda)$
coincide, so we may say that a hyperplane is defined by a point 
$a\in\TP^{l-1}$. The hyperplane defined by 
$(0,\dots,0)$ is the set of points $v\in \TP^{l-1}=\reals^{l-1}$ such that
$v$ either lies in the boundary of the positive orthant or has minimum
coordinate negative and repeated at least twice. Said in a more
compact (and invariant) form, it
equals the $l-2$-skeleton of the normal fan of the simplex with
vertices $O,e_1,\dots,e_{l-1}$. The translation of this hyperplane 
by the vector $-a\in\reals^l$ gives the hyperplane defined by $a$.

If $H$ is the hyperplane defined by a point $v\in
\TP^{l-1}$, here we call $-H$ the {\em anti-hyperplane} defined by $v$.
For the purposes of this paper,
the following consequence of the results in~\cite[Section 3]{DevStu}
can be taken as a definition of tropical
convex hull:

\begin{proposition}
\label{prop:tropical}
Let $v_1,\dots,v_k$ be a finite set of points in tropical $l-1$
space. Then, its tropical convex hull $\tconv(v_1,\dots,v_n)$ equals
the union of all bounded cells in the polyhedral arrangement of
tropical anti-hyperplanes given by $v_1,\dots,v_k$.
\end{proposition}

\begin{figure}[htb]
\begin{center}
\includegraphics[height=1in]{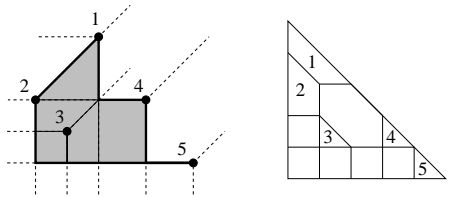}
\caption{
A tropical point configuration, and its associated regular mixed
subdivision of $k\Delta^{l-1}$}
\label{fig:tropical}
\end{center}
\end{figure}

The left part of 
Figure~\ref{fig:tropical} shows an example of this. The tropical
convex hull of the five dots equals the shaded region, including 
its boundary and the horizontal segment that reaches to point number 5.
Develin and Sturmfels make no clear distinction between the tropical convex
hull as a subset of $\TP^{l-1}$ and the polyhedral complex
in the above statement, and use the term ``tropical polytope'' referring to
both. Here we will use ``tropical polytope'' referring to the region and
call ``tropical order type'' 
of the point set the polyhedral complex. Two
point sets are {\em combinatorially equivalent} if they  have the
same bounded complex, in a labeled sense (with the label of each cell
indicating its relative position in each of the anti-hyperplanes. This
is essentially what Develin and Sturmfels call the ``type'' of a cell).

The connection to mixed subdivisions is given in the following
statement, paraphrased from Section 4 of~\cite{DevStu}:

\begin{theorem}
Let $M$ and $M'$ be two $k\times l$ real matrices. Then,
the columns of $M$ and $M'$ produce the same tropical order
type if and only if they produce the same regular mixed subdivision of
$k\Delta^{l-1}$ (where the $i$th column specifies the heights to
lift the vertices of the $i$th copy of $\Delta^{l-1}$).
\end{theorem}

\begin{corollary}
\label{coro:tropical}
There is the same number of order types of $k$ points in
tropical $l-1$ space as coherent mixed subdivisions of $k\Delta^{l-1}$;
that is, regular subdivisions of $\Delta^{k-1}\times \Delta^{l-1}$.
\end{corollary}

Figure~\ref{fig:tropical} illustrates the correspondence between a
tropical point configuration with 5 points in 2-space
and a mixed subdivision of $5\Delta^2$. As is easy to check, there
is a 1-to-1 dimension (and order) reversing correspondence 
between the cells defined by the tropical point set and the cells in
the mixed subdivision. Unbounded cells in the tropical point set
correspond to boundary cells in the mixed subdivision. Corresponding
cells are orthogonal to one another. Actually, the link of a cell in
the tropical point set is the normal fan of the corresponding cell in
the mixed subdivision.

\begin{figure}[htb]
\begin{center}
\includegraphics[height=6in]{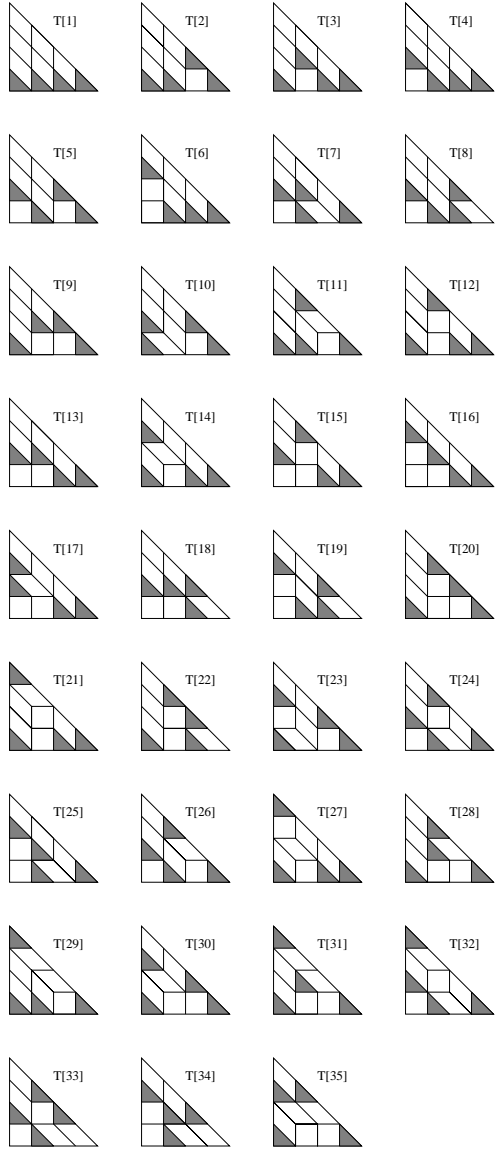}
\end{center}
\caption{The 35 symmetry classes of triangulations of $\Delta^2\times\Delta^3$}
\label{fig:35classes}
\end{figure}

As another  example, Figure~\ref{fig:35classes} 
shows all the fine mixed subdivisions of $3\Delta^2$, placed and
numbered to exactly match the 35 types of ``tropical quadrangles'' as
shown in Figure 6 of~\cite{DevStu}. (This list was originally computed by J.
Rambau~\cite{Rambau-TOPCOM}, and its order is the one given as 
output by TOPCOM,
that performs a ``breadth-first search'' on the graph of flips).

We now use
Corollary~\ref{coro:tropical} to give a bound on the number of regular
subdivisions of $\Delta^{k-1}\times \Delta^{l-1}$. For this, we extend
the tropical arrangement of $k$ anti-hyperplanes that defines
$\tconv(v_1,\dots,v_n)$ to a (usual) affine
arrangement $\H$ of $k{l \choose 2}$  hyperplanes. Indeed, the anti-hyperplane corresponding
to a point $v_i=(v_{i,1},\dots,v_{i,l-1})\in\reals^{l-1}\cong\TP^{l-1}$ is
a $(l-2)$-dimensional polyhedral complex with ${l\choose 2}$ maximal
cells, lying respectively in the following hyperplanes:
\begin{equation}
\begin{array}{l}
H(i;j,l):=\{x_j=v_{i,j}\},  \qquad \qquad \qquad  \text{for }  j=1,\dots,l-1, \text{ and }\\
H(i;j,k):=\{x_j-x_k=v_{i,j}-v_{i,k}\},   \qquad \text{for } 1\le j<k\le l-1. \\
\end{array}
\label{eqn:hyperplanes}
\end{equation}
Clearly,  point sets
with different tropical order type produce different (labeled)
hyperplane arrangements. Then:

\begin{theorem}
\label{thm:regular}
For any $k$ and $l$ the number of regular subdivisions of
$\Delta^{l-1}\times\Delta^{k-1}$ is bounded above by:
\[
\left(\frac{{\bf e}}{2}kl\right)^{l(l-1)(k-1)}
\]
\end{theorem}

\begin{proof}
We can assume that $(k-1)(l-1)\ge 6$, because if this is not the case, the bound can be proved by direct inspection:
if $l=2$ then the number of triangulations is $k!$ and if $l=k=3$ then the number of triangulations is 108
by Table \ref{table:numbers}.

We need to bound the number of different arrangements $\H$ that can be produced for varying 
$(v_1,\dots,v_k)\in \reals^{kl}$. Two arrangements are ``equal'' if they have the same chirotope,
that is to say, if every determinant of $l$ of the $k{l\choose 2}$ hyperplanes has the same sign in the two arrangements.

From the definition of the hyperplanes in equation (\ref{eqn:hyperplanes}) it is clear that
each of the determinants that define the chirotope of $\H$ is a linear functional on the $kl$ variables
$(v_{i,j})$. Since the tropical order type is invariant under addition of a constant to a row or column of the matrix 
$(v_{i,j})$, we can assume $v_{0,j}=v_{i,0}=0$ for every $i$ and $j$, leaving only $(l-1)(k-1)$ variables.




Hence, the order type of $\H$ appears represented as a cell in a 
huge linear hyperplane arrangement of ${k{l\choose 2} \choose l}\le \frac{k^l l^{2l}}{2^l l!}\le \left(\frac{{\bf e}}{2}kl\right)^l $ 
hyperplanes in $\reals^{(k-1)(l-1)}$. This gives the statement, since the number of cells in an 
arrangement of $N$ hyperplanes in $\reals^D$ is maximal for simple arrangements,
in which case it equals
\[
\sum_{i=0}^D 2^i {N\choose i} \le (D+1) 2^D {N\choose D} \le \frac{(D+1)2^D}{D!} N^D \le N^D.
\]
In the first inequality we assume that $D\le N/2$, which always happens for $N=\left(\frac{k l{\bf e}}{2}\right)^l$ and $D=(k-1)(l-1)$; 
in the last inequality we have used our assumption that $D=(k-1)(l-1)\ge 6$.
\end{proof}

Observe that our bound is quite rough
not only because different arrangements $\H$ may represent the same tropical order type, 
but also because the sign of many of the ${k{l\choose 2} \choose l}$ determinants considered in the proof is constant (independent of
the ${v_{i,j}}'s$). We believe the actual number of regular subdivisions to be in $(kl)^{\Theta(kl)}$.
Anyway, for fixed $l$ our bound gives the exact asymptotic behavior of the number of regular subdivisions:

\begin{corollary}
\label{coro:nonregular2}
For any fixed $l\ge 2$, the number of regular subdivisions of 
$\Delta^{l-1}\times \Delta^{k-1}$ is in
$
k^{\Theta(k)}.
$ For $l\ge  3$
the number of all subdivisions is in
$
2^{\Omega(k^2)}.
$
\end{corollary}

\begin{proof}
For regular subdivisions, the  upper bound follows from the previous theorem and the lower bound is trivial, since a single orbit of regular triangulations has already $k! \in k^{\Omega(k)}$ elements.

For all subdivisions, Theorem~\ref{thm:asymptotics} gives the case
$l=3$ and the others follow immediately: any subdivision of a
particular $\Delta^{2}\times \Delta^{k-1}$ face of $\Delta^{l-1}\times
\Delta^{k-1}$ can be extended to the whole polytope.
\end{proof}

\bigskip

\subsection*{Bibliographic remark:}
The writing of this paper has spanned an unusually long period of time, the first
drafts dating back to 1998. Previous versions of it have been cited as ``in preparation''
under the title {\em Applications of the polyhedral Cayley Trick to triangulations of polytopes}.

\bibliographystyle{amsplain}

\end{document}